\documentclass[leqno,11pt,twoside]{article}

\usepackage{amsfonts}                                %packages
\usepackage{amsmath}
\usepackage{amsthm}
\usepackage{amssymb}
\usepackage{latexsym}

\newtheorem{The}{Theorem}[section]              %theorem environments
\newtheorem{Pro}[The]{Proposition}
\newtheorem{Lem}[The]{Lemma}
\newtheorem{Cor}[The]{Corollary}

\theoremstyle{definition}
\newtheorem{Def}[The]{Definition}

\numberwithin{equation}{section}

\newcommand{\ra}{\rightarrow}
\newcommand{\lra}{\longrightarrow}
\newcommand{\sst}{\subseteq}
\newcommand{\mt}{\mapsto}

\newcommand{\ca}{\ensuremath{\mathcal{A}}}
\newcommand{\ce}{\ensuremath{\mathcal{E}}}

\newcommand{\cS}{\ensuremath{\mathcal{S}}}

\newcommand{\cu}{\ensuremath{\mathcal{U}}}
\newcommand{\Gak}{\ensuremath{G_{\ca}(k,n)}}
\newcommand{\Gan}{\ensuremath{G_{\ca}(n)}}
\newcommand{\gak}{\ensuremath{\mathcal{G}_{\ca}(k,n)}}
\newcommand{\gan}{\ensuremath{\mathcal{G}_{\ca}(n)}}
\newcommand{\Vak}{\ensuremath{V_{\ca}(k,n)}}

\newcommand{\vS}[2]{\ensuremath{\mathcal{V}^#1(#2)}}
\newcommand{\Nx}[1]{\ensuremath{\mathcal{N}(#1)}}

\newcommand{\C}{\mbox{$\mathbb{C}$}}
\newcommand{\K}{\mbox{$\mathbb{K}$}}
\newcommand{\N}{\mbox{$\mathbb{N}$}}
\newcommand{\R}{\mbox{$\mathbb{R}$}}

\newcommand{\tX}{\mbox{$\mathfrak{T}_X$}}

\newcommand{\gra}{\alpha}
\newcommand{\gl}{\lambda}
\newcommand{\gr}{\rho}
\newcommand{\Ua}{U_{\gra}}

\begin{document}

\renewcommand{\thefootnote}{\fnsymbol{footnote}}

 \pagestyle{myheadings}
\markboth{\sc M. H. Papatriantafillou --%
E.~Vassiliou}{\sc Classification of Vector Sheaves}

\title{\bf Grassmann sheaves and the classification of vector sheaves}

\author{\bf M.~H.~Papatriantafillou -- E.~Vassiliou\footnote{The authors
were partially supported by University of Athens Research Grands
70/4/5639 and 70/4/3410, respectively.}}

\date{}
\maketitle \thispagestyle{empty}

\begin{abstract}\noindent
Given a sheaf of unital commutative and associative algebras \ca,
first we construct the $k$-th Grassmann sheaf $\gak$ of $\ca^n$,
whose sections induce vector subsheaves of $\ca^n$ of rank $k$. Next
we show that every vector sheaf over a paracompact space is a
subsheaf of $\ca^{\infty}$. Finally, applying the preceding results
to the universal Grassmann sheaf $\gan$, we prove that vector
sheaves of rank $n$ over a paracompact space  are classified by the
global sections of $\gan$.
\end{abstract}

%--------------------------------------------------------------------
%\maketitle
%--------------------------------------------------------------------
%%%%%%%%%%%%%%%%%%%%%%%%%%%%%%%%%%%%%%%%%%%%%%%%%%%%%%%%%%%%%%%%%%%%%

\section*{Introduction}                        %Introduction

\noindent Let \ca\ be a sheaf of unital commutative and associative
algebras over the ring \R\ or \C. A vector sheaf \ce\ is a locally
free \ca-module. For instance, the sections of a vector bundle
provide such a sheaf. However, a vector sheaf is not necessarily
free, as is the case of the sections of a non trivial vector bundle.

Recently, vector sheaves gained a particular interest because they
serve as the platform to abstract the classical geometry of vector
bundles and their connections within a non smooth framework. This
point of view has already been developed in \cite{Mal1} (see also
\cite{Mal2} for applications to physics, and \cite{Vas} for the
reduction of the geometry of vector sheaves to the general setting of
principal sheaves).

A fundamental result of the classical theory is the homotopy
classification of vector bundles (of rank, say, $n$) over a fixed
base. The construction of the classifying space, and the subsequent
classification, are based on the Grassmann manifold (or variety)
$G_k(\R^n)$ of $k$-dimensional subspaces of $\R^n$. In this respect
we refer, e.g., to \cite{Hus} and \cite{Kah}. However, considering
vector sheaves, we see that a homotopy classification is not
possible, since the pull backs of a vector sheaf by homotopic maps
need not be isomorphic, even in the trivial case of the free
$\ca$-module $\ca$, as we prove in Section 1. Consequently, any
attempt to classify vector sheaves (over a fixed space $X$) should
not involve pull-backs and homotopy.

In this paper we develop a classification scheme based on a sort of
universal Grassmann sheaf. More explicitly, for fixed $k \leq n \in
\N$, in Section 2 we construct --in two equivalent ways--  a sheaf
\gak, legitimately called the $k$-th Grassmann sheaf of $\ca^n$,
whose sections coincide (up to isomorphism) with vector subsheaves
of $\ca^n$ of rank $k$ (Proposition 2.3). Then, inducing in Section
3 the vector sheaf $\ca^{\infty}$, we show that every vector sheaf
over a paracompact space is a subsheaf of $\ca^{\infty}$ (Theorem
3.1). A direct application of the previous ideas leads us to the
construction of the universal Grassmann sheaf $\gan$ of rank $n$.
The main result here (Theorem 3.5) asserts that arbitrary vector
sheaves of rank $n$, over a paracompact base space, coincide --up to
isomorphism-- with the sections of $\gan$.

\section{Vector sheaves and homotopy}                %Section 1

\noindent For the general theory of sheaves we refer to standard
sources such as \cite{Bre}, \cite{Dow}, \cite{God}, and \cite{Ten}.
In what follows we recall a few definitions in order to fix the
notations and terminology of the present paper.

\smallskip  Throughout the paper \ca\ denotes a fixed sheaf of {\em unital
commutative and associative \K-algebras} ($\K = \R, \C$) over a
topological space $X$. An \ca-\emph{module} $\ce \equiv (\ce,\pi,X)$
is a sheaf whose stalks $\ce_x$ are $\ca_x$-modules so that the
respective operations of addition and scalar multiplication
\[
   \ce\times_X \ce \lra  \ce \quad \text{and}
   \quad \ca\times_X \ce \lra  \ce
\]
are continuous. In particular, a \emph{vector sheaf of rank n} is an
\ca-module \ce, locally isomorphic to $\ca^n$. This means there is an
open covering $\cu = \{\Ua\},\; \gra\in I$, of $X$ and
$\ca|_{\Ua}$-iso\-mor\-phisms
\[
    \psi_{\gra}\colon \ce|_{\Ua} \lra \ca^n|_{\Ua},
                                              \qquad \gra\in I.
\]
The category of vector sheaves of rank $n$ over $X$ is denoted by
\vS{n}{X}. More details, examples and applications of vector sheaves
can be found in \cite{Mal1}, \cite{Mal2}, and \cite{Vas}.

As already mentioned in the Introduction, we shall show, by a
concrete counterexample, that homotopic maps do not yield isomorphic
pull-backs, even in the simplest case of the free \ca-module \ca. In
fact, we consider two {\em non-isomorphic} algebras $A_0$ and $A_1$
and a morphism of algebras $\gr \colon A_0 \ra A_1$. Given now a
topological space $X$ and a fixed point $x_0\in X$, for every open $U
\sst X$ we set
\[
    A(U):=
        \begin{cases}
      A_0, & x_0 \in U,\\
      A_1, & x_0 \notin U,
        \end{cases}
\]
while, for every open $V \sst U$, $\gr^U_V \colon A(U) \ra A(V)$
denotes the corresponding (restriction) map, defined  by
\[
    \gr^U_V :=
    \begin{cases}
      id : A_0 \ra A_0,  & \text{if}\quad x_0 \in V,\\
     \,\,\gr : A_0 \ra A_1, & \text{if}\quad U\ni x_0 \notin V,\\
      id : A_1 \ra A_1,  & \text{if}\quad x_0 \notin U.
 \end{cases}
\]

It is not difficult to show that $\big(A(U),\gr^U_V\big)$ is a
presheaf whose sheafification is a sheaf of algebras, denoted by $\ca
\equiv (\ca,\pi,X)$. It is clear that
\[
     \ca_{x_0} = \lim_{\stackrel{\displaystyle{\lra}}{U\in
                             \Nx{x_0}}}   A(U) = A_0.
\]
On the other hand, if there is an $x_1\in X$ admitting a neighborhood
$V\in \Nx{x_1}$ with $x_0 \notin V$ (a fact always ensured if $X$ is
a $T_1$-space), then
\[
    \ca_{x_1} = \lim_{\stackrel{\displaystyle{\lra}}{x_1 \in
                  W \sst V}}  A(W) = A_1,
\]
which is not isomorphic to $\ca_{x_0}$. Hence, a vector sheaf, even a
free one, need not have locally isomorphic fibres.

Let now $\gra \colon [0, 1] \ra X$ be a continuous path with $\gra(0)
= x_0$ and $\gra(1) = x_1$. For any topological space $Y$, we define
the map
\[
   f\colon [0, 1] \times Y  \lra X \quad \text{with} \quad
    f(t,y) := \gra(t).
\]
Obviously, this is a homotopy between the constant maps $f_0 =x_0
\colon Y\ra X$ and $f_1=x_1\colon Y\ra X$. As a result, taking the
pull-backs of \ca\ by the latter, we see that, for every $y\in Y$,
\begin{align*}
f_0^{\,*}(\ca)_y &= \{(y,a)\,|\, a\in \ca_{x_0} \} = \{y\}\times A_0,\\
f_1^{\,*}(\ca)_y &= \{(y,b)\,|\, b\in \ca_{x_1} \} = \{y\}\times A_1;
\end{align*}
that is, we obtain two non-isomorphic stalks, thus proving the claim.

\section{The Grassmann sheaf of rank $k$ in $\ca^n$}       %Section 2

\noindent As in Section 1, $\ca$ is a sheaf of unital commutative
and associative \K-algebras over a given topological space $X$. We
denote by \tX\ the topology of $X$ and fix $n\in \N$. For $k\in \N$
with $k\leq n$ and any $U\in \tX$, we define the set
\[
   \Gak(U) := \big\{\cS\; \text{subsheaf of}\; \ca^n|_U \,:\,
                                        \cS \cong \ca^k|_U \big\},
\]
that is, $\Gak(U)$ consists of the free submodules of $\ca^n|_U$ of
rank $k$. If, for every $U,V\in \tX$ with $V\sst U$,
\[
   \gr^U_V \colon \Gak(U)\lra \Gak(V)\colon \cS\mapsto \cS|_V
\]
denotes the natural restriction, it is clear that the collection
\begin{equation}                                             %(2.1)
    \Gak := \left(\Gak(U),\, \gr^U_V\right)
\end{equation}
determines a presheaf. Moreover, it is a {\em monopresheaf}. Indeed,
if $U = \bigcup _{i\in I} U_i$ and $\ce_1, \ce_2 \in
G_{\ca}(n,k)(U)$ with $\ce_1|_{U_i} = \ce_2|_{U_i}$, for all $ i \in
I$, then $\ce_1 = \ce_2$. However, it is not complete: If $\ce_i \in
G_{\ca}(n,k)(U_i)$, with $\ce_i|_{U_i\cap U_j} = \ce_j|_{U_i\cap
U_j}$, then $\ce := \bigcup_{i\in I} \ce_i$ is a vector sheaf over
$U$, but not necessarily a free $\ca|_U$-module.

\begin{Def}                                           %Definition 2.1
The \emph{$k$-th Grassmann sheaf of $\ca^n$}, denoted by \gak, is
defined to be the sheaf generated by the presheaf $\Gak$.
\end{Def}

Since $\Gak$ is not complete, it does not coincide with the complete
presheaf
 \begin{equation}                                            %(2.2)
 \left(\gak(U), r^U_V\right),
 \end{equation}
of (continuous) sections of \gak. We shall describe $\gak$ via
another complete presheaf. As a matter of fact, we consider the
presheaf
\begin{equation}                                             %(2.3)
    \Vak := \big(\Vak(U),\, \gl^U_V\big),
\end{equation}
where now
\[
   \Vak(U) := \big\{\cS\; \text{subsheaf of}\; \ca^n|_U\,:\,
                                         \cS\in \vS{k}{U}\big\}
\]
and
\[
   \gl^U_V \colon \Vak(U)\lra \Vak(V)\colon \cS\mapsto \cS|_V
\]
are the natural restrictions. In contrast to \Gak, \Vak\ is obviously
a complete presheaf.

\begin{Lem}
The sheaf generated by $\Vak$ is isomorphic to $\gak$.
\end{Lem}

\begin{proof}
Clearly, for every $U \in \tX$,
\[
\Gak(U) \sst \Vak(U),
\]
that is,\ \Gak \ is a sub-presheaf of \Vak. Besides, for every $\ce
\in \Vak(U)$ and every $x \in U$, there is $V \in \tX$ with $x \in V
\sst U$, so that $\ce|_V$ is free of rank $k$, namely $\ce|_V \in
\Gak(V)$. Thus, $\Gak(U)$ and $\Vak(U)$ define the same sheaf $\gak$.
\end{proof}

Since $\Vak$ is complete, it is isomorphic with the sheaf of sections
of $\gak$, thus we have the following interpretation of the elements
of $\gak(X)$.

\begin{Pro}
The global sections of the $k$-th Grassmann sheaf $\gak$ coincide
--up to isomorphism-- with the vector subsheaves of $\ca^n$ of rank
$k$.
\end{Pro}

\section{The universal Grassmann sheaf}                  %Section 3

\noindent The preliminary results of the preceding section hold for
every base space $X$. Here we prove that if $X$ is a paracompact
space, then {\em any} vector sheaf can be interpreted as a section
of an appropriate universal Grassmann sheaf.

\smallskip
First we prove a Whitney-type embedding theorem. To this end, for
every sheaf of algebras $\ca$, we consider the presheaf
 \[
 U \longmapsto \prod_{i\in \N} \ca_i(U), \qquad U \in \tX,
 \]
where $\ca_i = \ca$, for every $i \in \N$, with the obvious
restrictions. This presheaf generates the infinite {\em fibre}
product
 \[
 \ca^{\infty} := \prod_{i\in \N} \ca_i,
 \]
which is a free $\ca$-module. Then, we obtain:

\begin{The}                                             %Theorem 3.1
Let $X$ be a paracompact space. Then every vector sheaf $\ce$ of
finite rank over $X$ is a subsheaf of $\ca^{\infty}$.
\end{The}

\begin{proof}
Let $\ce$ be a vector sheaf of rank, say, $k$. Since $X$ is
paracompact, a reasoning similar to that of \cite[Proposition
5.4]{Hus} proves that $\ce$ is free over a countable open covering
$\{U_i\}_{i \in \N}$ of $X$. Let $\psi_i \colon \ce|_{U_i} \ra
\ca^k|_{U_i}$, $i \in \N$, be the respective family of $\ca$-module
isomorphisms. The same open covering has a countable locally finite
open refinement (\cite[Ch.\ VIII, Theorem 1.4]{Dug}), with a
subordinate partition of unity $\{\alpha_i \colon X \ra \R\}_{i \in
\N}$ (ibid., Ch.\ VIII, Theorem 4.2). For every $i \in \N$, we
define the map $\alpha_i\psi_i \colon \ce \ra \ca^k$ by
 \[
 \alpha_i\psi_i(u) := \begin{cases}
                      \alpha_i(\pi(u))\psi_i(u), & \pi(u) \in U_i,\\
                       0, & \pi(u) \notin U_i.
                       \end{cases}
                      \]
Therefore, $\alpha_i\psi_i$ is an $\ca$-module morphism, whose
restriction to the interior of ${\rm supp} \, \alpha_i \sst U_i$ is
an isomorphism.

We consider the fibre product $\prod_{i\in \N}(\ca^k)_i$, where
$(\ca^k)_i \equiv \ca^k$, for every $i \in \N$, and we denote by
 \[
 p_i \colon \prod_{i\in \N}(\ca^k)_i \lra (\ca^k)_i
 \]
the corresponding projections. The universal property of the product
ensures the existence of a unique $\ca$-morphism
 \[
 \psi : \ce \lra \prod_{i\in \N}(\ca^k)_i,
 \]
such that
 \[
 p_i \circ \psi = \alpha_i\psi_i.
 \]
Then $\psi$ is a monomorphism. In fact, let $0 \neq u \in \ce_x$ with
$\psi_x(u) = 0$. There is $i \in \N$, with $\alpha_i(x)
> 0$, thus $\alpha_i(x)\psi_{i,x}(u) \neq 0$, a
contradiction. Hence, $\ce$ is identified with its image $\psi(\ce)
\leq \prod_{i\in \N}(\ca^k)_i$. Since
 \[
 \prod_{i\in \N}(\ca^k)_i \equiv \prod_{i\in \N} \ca_i,
 \]
where $\ca_i = \ca$, for every $i \in \N$, the assertion is proven.
\end{proof}

We shall show that a further restriction on the topology of $X$
leads to an embedding of $\ce$ into a smaller sheaf. To this end,
assume that $\ce$ is a vector sheaf of rank $k$, which is free over
a {\em finite} open covering $\{ U_i \}_{1\leq i\leq n}$ of $X$. Let
$\psi_i \colon \ce|_{U_i} \ra \ca^k|_{U_i}$ be the respective
$\ca$-module isomorphisms, and $\{\alpha_i \colon X \ra \R\}_{1 \leq
i \leq n}$ a subordinate partition of unity. Considering the maps
$\alpha_i\psi_i \colon \ce \ra \ca^k$, as before, we obtain the
sheaf morphism
 \[
 f \colon \ce \lra \ca^{kn} : u \mt \big(\alpha_1(\pi(u))\cdot\psi_1(u),
 \dots,\alpha_n(\pi(u))\cdot\psi_n(u)\big)
 \]
which embeds $\ce$ into the free $\ca$-module $\ca^{kn}$. Therefore
we have proved the following:

\begin{Pro}
Over a compact space $X$, every vector sheaf of finite rank is a
subsheaf of a free $\ca$-module of finite rank.
\end{Pro}

Clearly, every free $\ca$-module $\ca^n$ is a submodule of the free
$\ca$-module $\bigoplus_{i \in \N} \ca_i$, with $\ca_i = \ca$, for
every $i \in \N$. Thus we obtain:

\begin{Cor}
Over a compact space $X$, every vector sheaf of finite rank is a
subsheaf of the free $\ca$-module $\bigoplus_{i \in \N} \ca_i$.
\end{Cor}

We are now in a position to repeat the constructions of Section 2 in
a more general way. For every $n\in \N$, we define the set
 \[
 \Gan(U) := \big\{\cS\; \text{subsheaf of}\; \ca^{\infty}|_U \,:\,
                                        \cS \cong \ca^n|_U \big\}.
 \]
Then the collection
\[
    \Gan := \big(\Gan(U),\, \gr^U_V\big), \qquad U\in\tX,
\]
where $\gr^U_V$ denotes the obvious restriction, is a non-complete
monopresheaf.

\begin{Def}
The sheaf $\gan$, generated by the presheaf $\Gan$, is called {\em
the universal Grassmann sheaf of rank $n$}.
\end{Def}

The respective complete presheaf of the sections of $\gan$ is
isomorphic to the presheaf
\[
    V_{\ca}(n) := \big(V_{\ca}(n)(U),\, \gl^U_V\big), \qquad U\in\tX,
\]
where now
\[
   V_{\ca}(n)(U) := \big\{\cS\; \text{subsheaf of}\; \ca^{\infty}|_U\,:\,
                                         \cS\in \vS{n}{U}\big\}
\]
and $\gl^U_V$ are the natural restrictions.

\medskip
As a result, adapting the proof of Proposition 2.3 to the present
situation, we obtain the main result of this work, namely the
following classification of vector sheaves:

\begin{The}
If $X$ is a paracompact space, then the vector sheaves of rank $n$
(over $X$) coincide --up to isomorphism-- with the global sections
of the universal Grassmann sheaf ${\mathcal G}_{\ca}(n)$.
\end{The}

%%%%%%%%%%%%%%%%%%%%%%%%%%%%%%%%%%%%%%%%%%%%%%%%%%%%%%%%%%%%%%%%%%%%%

\end{document}